\documentclass[11pt]{article}
\usepackage{graphicx}
\usepackage{amsmath}
\usepackage[dvips]{epsfig}
\usepackage{amssymb}
\begin{document}
\begin{center} {\Large \bf  On the  $(h,q)$-zeta function associated with  $(h,
q)$-Bernoulli numbers and polynomials}
\\ \vspace*{12 true pt}  Taekyun   Kim
\vspace*{12 true pt} \\
 Division of General Education-Mathematics,\\
 Kwangwoon University, Seoul 139-701,  Korea \\
               e-mail: tkkim@kw.ac.kr
          \end{center}
%--------------------------------------------------------------------
\vspace*{12 true pt} \noindent {\bf Abstract :} In this paper we
consider the following $(h, q)$-zeta function:
$$ \zeta_q(s,x \mid h)= \sum_{n=0}^\infty \dfrac{q^{hn+x}} { [n+x]_q^s}
+ \dfrac{(q-1)(h-1)}{1-s} \sum_{n=0}^\infty \dfrac{q^{(h-1)n}}{
[n+x]_q^{s-1}}, $$ where $ x \neq 0, -1, -2, \cdots,  s \in \Bbb C
\setminus \{ 1 \} $ and  $ h \in \mathbb{C}. $ Finally, we lead to
a useful integral representation for the $(h, q)$-zeta functions
and give the functional equation associated with $(h,
q)$-Bernoulli numbers and polynomials.

\vspace*{12 true pt} \noindent {\bf 2000 Mathematics Subject
Classification :} 11B68, 11S40, 11S80

 \vspace*{12 true pt} \noindent {\bf
Key words :}   Bernoulli  numbers, Bernoulli  polynomials,
 $(h,q)$-Bernoulli numbers and polynomials,  $(h,q)$-zeta function

\begin{center} {\bf 1. Introduction } \end{center}

When one talks of $q$-extension, $q$ is considered  in many ways
such as an indeterminate, a complex number $q\in  \mathbb{C},$ or
$p$-adic number $q\in \Bbb C_p .$ Throughout this paper we assume
that $q\in \Bbb C$ with $|q|<1$.
 We use the notation  of $q$-number as $$[x]_q =\frac{1-q^x}{1-q}.$$
Note that $\lim_{ q \rightarrow 1} [x]_q =x,$( see [6]).

The Bernoulli polynomials in $\Bbb C$ are defined by the formula
$$ e^{xt} \dfrac{t}{e^t-1}= e^{B(x)t}= \sum_{n=0}^\infty B_n (x) \dfrac{t^n}{n!},$$
with the usual convention of replacing $B^i(x)$ by $B_i(x)$, (see
[1-14]). In the special case, $x=0, B_n(0)=B_n$ are called the
$n$-th ordinary Bernoulli numbers, (see [1-10]). The Bernoulli
numbers are used to express special values of Riemann zeta
function, which is defined by
$$ \zeta(s)= \sum_{n=1}^\infty  \dfrac{1}{n^s}, s \in \Bbb C,
(\text{ see} [7]).$$ That is, $\zeta(2m)=
\dfrac{(2\pi)^{2m}(-1)^{m-1} B_{2m}}{2( 2m)!},$ for $m \in \Bbb
N,$ and $\zeta(1-2m)= -\dfrac{B_{2m}}{2m},$ (see [7-14]). First,
we consider $(h, q)$-Bernoulli numbers and polynomials as the
$q$-extension of Bernoulli numbers and polynomials. From those
numbers and polynomials, we derive some $(h, q)$-zeta functions as
the $q$-extension of Riemann zeta function. That is, the purpose
of this paper ir to study the following $(h, q)$-zeta function:
$$ \zeta_q(s,x \mid h)= \sum_{n=0}^\infty \dfrac{q^{hn+x}} { [n+x]_q^s}
+ \dfrac{(q-1)(h-1)}{1-s} \sum_{n=0}^\infty \dfrac{q^{(h-1)n}}{
[n+x]_q^{s-1}}, $$ where $ x \neq 0, -1, -2, \cdots,$ and $  h \in
\mathbb{C}, $   $  s \in \Bbb C \setminus \{ 1 \} . $

Finally, we derive to a useful integral representation for the
$(h, q)$-zeta functions and give the functional  equation
associated with $(h, q)$-Bernoulli numbers and polynomials.
Recently,  several authors have studied the $q$-zeta functions and
the $q$-Bernoulli numbers( see [1-14]). Our $q$-extensions of
Bernoulli numbers and polynomials in this paper are different the
$q$-extension of  Bernoulli numbers and polynomials  which are
treated by several authors in previous papers.

\bigskip
\begin{center} {\bf 2. $(h, q)$-Bernoulli numbers and polynomials associated with  $(h, q)$-zeta functions} \end{center}
\bigskip

For $ h \in \mathbb{C}, $ let us consider $(h, q)$-Bernoulli
polynomials as follows:

$$  \aligned     F_q(t, x \mid h) & = -t \sum_{m=0}^\infty q^{hm +x} e^{[x+m]_q t}
+ (h-1)(1-q) \sum_{m=0}^\infty q^{(h-1)m} e^{[x+m]_q t}  \\
&  = \sum_{n=0}^\infty \beta_{n, q}^h(x) \dfrac{t^n}{n!}.
\endaligned  \eqno(1)$$
From (1), we note that

$$ F_q(t, x \mid h) = e^{\left( \dfrac{1}{1-q} \right) t} \sum_{l=0}^\infty (-1)^l q^{lx} \dfrac{l+h-1}{[l+h-1]_q}
\left( \dfrac{1}{1-q} \right)^l \dfrac{t^l}{l!}. \eqno(2)$$

By (1) and (2), we get
$$\beta_{n, q}^h(x)=\dfrac{1}{(1-q)^n}  \sum_{l=0}^n  \binom nl (-1)^l q^{lx} \dfrac{l+h-1}{[l+h-1]_q}. \eqno(3)$$
In the special case,  $x=0$, $\beta_{n, q}^h (0)= \beta_{n,  q
}^h$ are called the $n$-th  $(h, q)$-Bernoulli numbers.

In (3), it is easy to show that
$$\beta_{n, q}^h(x)=( [x]_q + q^x \beta ^h )^n, \eqno(4) $$
where we use the usual convention about replacing $
(\beta^h)^n=\beta_{n, q}^h.$

Note that
$$\lim_{ q \rightarrow 1} F_q(t, x \mid h) = \dfrac{ e^{xt} t}{e^t-1}= \sum_{n=0}^\infty  B_n (x) \dfrac{t^n}{n!}.$$
Let $      F_q(t  \mid h)  = \sum_{n=0}^\infty \beta_{n, q}^h
\dfrac{t^n}{n!}$ be the generating function of $(h, q)$-Bernoulli
numbers. Then we get
$$   F_q(t  \mid h)  = -t \sum_{m=0}^\infty q^{hm } e^{[m]_q t}
+ (h-1)(1-q) \sum_{m=0}^\infty q^{(h-1)m} e^{[m]_q t} . \eqno(5)$$
From (1), (4), and (5), we can derive the following difference
equation:

$$  \aligned     F_q(t, x \mid h) & =  e^{[x]_q t} F_q( q^x t  \mid h) \\
& = -t \sum_{n=0}^\infty q^{hn +x} e^{[x+n]_q t} + (h-1)(1-q)
\sum_{n=0}^\infty q^{(h-1)n} e^{[x+n]_q t}.
\endaligned  $$
Therefore, we obtain the following proposition.

\medskip
{ \bf Proposition 1.} For $ h  \in \Bbb C,$  we have
$$ \aligned     \beta_{n, q}^h(x)  & =  \dfrac{1}{(1-q)^n}  \sum_{l=0}^n \binom nl (-1)^l q^{lx} \dfrac{l+h-1}{[l+h-1]_q} \\
&=  \sum_{l=0}^n \binom nl  q^{lx}  \beta_{l, q}^h [x]_q^{n-l} .
\endaligned
 $$
\medskip

Note that
$$\left. \left( \frac{d}{dt}\right)^k F_q(t, x \mid h) \right\vert_{t=0}
= -k  \sum_{n=0}^\infty q^{hn +x} {[x+n]_q }^{k-1} + (h-1)(1-q)
\sum_{n=0}^\infty q^{(h-1)n} {[x+n]_q}^k.$$ Thus, we obtain the
following corollary.

\medskip
{ \bf Corollary 2.} For $ k  \in \Bbb Z_+= \Bbb N \cup \{ 0 \},$
we have
$$  \beta_{k, q}^h (x)  =
-k  \sum_{n=0}^\infty q^{hn +x} [x+n]_q^{k-1} + (h-1)(1-q)
\sum_{n=0}^\infty q^{(h-1)n} [x+n]_q^k ,
 $$
  and
$$  \beta_{k, q}^h   =
-k  \sum_{n=0}^\infty q^{hn } [n]_q^{k-1} + (h-1)(1-q)
\sum_{n=0}^\infty q^{(h-1)n} [n]_q^k . $$

\medskip

It is easy to show that
$$  q^{(h-1)n}  F_q(t, n \mid h)-   F_q(t \mid h) =  t \sum_{l=0}^{n-1} q^{hl}
e^{[l]_qt}-(h-1)(1-q)\sum_{l=0}^{n-1} q^{(h-1)l} e^{[l]_q t} .
 $$

Thus, we have
$$\beta_{0, q}^h =\dfrac{h-1}{[h-1]_q}, \text{ and }
 q^{h-1} \beta_{n, q}^h(1)- \beta_{n, q}^h =\delta_{1 n}, \eqno(6)$$
where $\delta_{1 n}$ is kronecker symbol.

Therefore, we obtain the following theorem.

\bigskip { \bf Theorem 3.} For $ h  \in \Bbb C, n \in \Bbb N$ and
$ m \in \Bbb Z_+$,  we have
$$ q^{(h-1)n} \beta_{m, q}^h(n) -\beta_{m, q}^h =
m  \sum_{l=0}^{n-1} q^{hl} [l]_q^{m-1}-(h-1)(1-q)\sum_{l=0}^{n-1}
q^{(h-1)l} [l]_q^m .
 $$
 In the special case, $n=1$, we have
$$\beta_{0, q}^h =\dfrac{h-1}{[h-1]_q}, \text{ and }
 q^{h-1} \beta_{n, q}^h(1)- \beta_{n, q}^h =\delta_{1 n}, $$
 where $\delta_{1 n}$ is kronecker symbol.

\bigskip

Now, we consider the following integral representation in complex
plane.
$$ \aligned   & \dfrac{1}{\Gamma(s)} \int_0^\infty t^{s-2} F_q(-t, x \mid h)dt \\
&= \sum_{n=0}^\infty \dfrac{ q^{hn +x}}{ [x+n]_q^s} +  \dfrac{\Gamma(s-1)}{\Gamma(s)}
(h-1)(1-q) \sum_{n=0}^\infty \dfrac{ q^{(h-1)n }}{ [x+n]_q^{s-1}}\\
&=  \sum_{n=0}^\infty \dfrac{ q^{hn +x}}{ [x+n]_q^s} +
\dfrac{(h-1)(1-q)}{s-1}  \sum_{n=0}^\infty \dfrac{ q^{(h-1)n }}{
[x+n]_q^{s-1}},
\endaligned
\eqno(7)
 $$
 where $ x \neq 0, -1, -2, \cdots,  h \in
\mathbb{C}, $  and  $  s \in \Bbb C \setminus \{ 1 \} . $

By the definition of $(h, q)$-Bernoulli polynomials, we see that
$$ \dfrac{1}{\Gamma(s)} \int_0^\infty t^{s-2} F_q(-t, x \mid h)dt
= \sum_{m=0}^\infty \left( \dfrac{1}{\Gamma(s)} \int_0^\infty
t^{m+s-2}dt \right) \dfrac{(-1)^m \beta_{m, q}^h(x)}{m!}.
\eqno(8)$$ In the special case, $s=1-k ( k \in \Bbb N)$, we see
form (7) and (8) and the basic theory of complex analysis
including Laurent series that
$$\sum_{n=0}^\infty  q^{hn +x} [x+n]_q^{k-1} -
\dfrac{(h-1)(1-q)}{k}  \sum_{n=0}^\infty  [x+n]_q^k q^{(h-1)n } =
- \dfrac{ \beta_{k, q}^h(x)}{k}. \eqno(9)$$
 In the viewpoint of
(7), we can define the following Hurwitz's type $(h, q)$-zeta
function:

\medskip
{ \bf Definition 4.}
 For $ h  \in \Bbb C,  s \in \Bbb C \setminus \{ 1 \}, $ and  $ x \neq 0, -1, -2, \cdots,  $
 define
$$ \zeta_q( s, x \mid h)=  \sum_{n=0}^\infty \dfrac{ q^{hn +x}}{ [x+n]_q^s} +
\dfrac{(h-1)(1-q)}{s-1}  \sum_{n=0}^\infty \dfrac{ q^{(h-1)n }}{
[x+n]_q^{s-1}}.
 $$
\medskip

Note that $\lim_{ q \rightarrow 1}\zeta_q( s, x \mid h)=\zeta( s,
x )$, where $ \zeta( s, x )= \sum_{n=0}^\infty \dfrac{1}{(n+x)^s}$
is called Hurwitz's zeta function.

\medskip
{ \bf Remark.}
 Note that $ \zeta_q( s, x \mid h)$ has only simple pole at $s=1$
 and $ \zeta_q( s, x \mid h)$ is meromorphic function except for
 $s=1$ in complex  $s$-plane.

\medskip
By (9) and Definition 4, we obtain the following theorem.

\bigskip { \bf Theorem 5.} For $ k \in \Bbb N$ ,   we have
$$ \zeta_q( 1-k, x \mid h)= - \dfrac{ \beta_{k, q}^h(x)}{k}.$$

\bigskip
In the special case, $x=1$, we see that
$$ \aligned    \zeta_q( s, 1 \mid h) &= \sum_{n=0}^\infty \dfrac{ q^{hn +1}}{ [n+1]_q^s} +
\dfrac{(h-1)(1-q)}{s-1}  \sum_{n=0}^\infty \dfrac{ q^{(h-1)n }}{
[n+1]_q^{s-1}}\\
&= \sum_{n=1}^\infty \dfrac{ q^{hn +1-h}}{ [n]_q^s} +
\dfrac{(h-1)(1-q)}{s-1}  \sum_{n=1}^\infty \dfrac{ q^{(h-1)n+1-h
}}{ [n]_q^{s-1}} \\
&= q^{-(h-1)} \left(  \sum_{n=1}^\infty \dfrac{ q^{hn}}{ [n]_q^s}
+ \dfrac{(h-1)(1-q)}{s-1}  \sum_{n=1}^\infty \dfrac{ q^{(h-1)n}}{
[n]_q^{s-1}}  \right).
\endaligned
 $$
Now, we define the $(h,q)$-zeta function as follows:

\medskip
{ \bf Definition 6.}
 For $s \in \Bbb C \setminus \{ 1 \},$ and $  h  \in \Bbb C$,
 define
$$ \zeta_q( s \mid h)=  \sum_{n=1}^\infty \dfrac{ q^{hn}}{ [n]_q^s} +
\dfrac{(h-1)(1-q)}{s-1}  \sum_{n=1}^\infty \dfrac{ q^{(h-1)n }}{
[n]_q^{s-1}}.
 $$
\medskip

Note that $$ \zeta_q( s \mid h)= q^{h-1} \zeta_q( s, 1  \mid h).$$
For $k \in \mathbb{N}$ with $k>1$, we have
$$ \zeta_q( 1-k \mid h)= q^{h-1} \zeta_q( 1-k, 1  \mid h)= - \dfrac{ q^{h-1} \beta_{k, q}^h(1)}{k}. \eqno(10)$$

By (6), (10), and Corollary 2, we obtain the following corollary.

\bigskip
 { \bf Corollary 7.} For $ k \in \Bbb N$ ,   we have
$$ \zeta_q( 1-k \mid h)= - \dfrac{ \beta_{k, q}^h}{k}.$$

\bigskip

Let $\chi$ be the  Dirichlet  character with conductor $f \in \Bbb
N$. Then  we define the generalized $(h,q)$-Bernoulli polynomials
attached to $\chi $ as follows:
$$\beta_{n,\chi,  q}^h(x)=  [f]_q^{n-1} \sum_{a=0}^{f-1} \chi(a) q^{(h-1)a} \beta_{n, q^f}^h \left( \dfrac{x+a}{f} \right).
\eqno(11)$$ Note that
$$\lim_{q \rightarrow 1} \beta_{n,\chi,  q}^h(x)=  f^{n-1} \sum_{a=0}^{f-1} \chi(a)  \beta_{n} \left( \dfrac{x+a}{f} \right)
=B_{n, \chi}(x),$$ where  $B_{n, \chi}(x)$ are the $n$-th
generalized ordinary Bernoulli polynomials attached to $\chi$.
 In the special
case,  $x=0$, $\beta_{n,\chi,  q}^h(0)= \beta_{n,\chi,  q}^h$ are
called the $n$-th generalized  $(h,q)$-Bernoulli numbers attached
to $\chi$.

From (3) and (11), we note that
$$ \beta_{n,\chi,  q}^h(x)=   \sum_{a=0}^{f-1} \chi(a) q^{(h-1)a}\dfrac{1}{(1-q)^n}
 \sum_{l=0}^n  \binom nl (-1)^l q^{l(x+a)} \dfrac{l+h-1}{[f(l+h-1)]_q}. \eqno(12)$$
By (12), we easily get
$$ \beta_{n,\chi,  q}^h(x)= (h-1)(1-q)
\sum_{m=0}^\infty q^{(h-1)m}\chi(m) [x+m]_q^k -n \sum_{m=0}^\infty
q^{hm +x} \chi(m) [x+m]_q^{n-1}. $$ Therefore, we obtain the
following theorem.

\bigskip
 { \bf Theorem 8.}  For $ h \in \Bbb C$  and $n \in \Bbb Z_+$ ,   we have
$$ \beta_{n,\chi,  q}^h(x)= (h-1)(1-q)
\sum_{m=0}^\infty q^{(h-1)m}\chi(m) [x+m]_q^k -n \sum_{m=0}^\infty
q^{hm +x} \chi(m) [x+m]_q^{n-1}. $$

\bigskip

 Let $      F_{\chi, q}(t, x  \mid h)  = \sum_{n=0}^\infty \beta_{n, \chi,
 q}^h(x)
\dfrac{t^n}{n!}.$  Then we see that
$$   F_{\chi, q}(t, x  \mid h) =(h-1)(1-q)
\sum_{m=0}^\infty q^{(h-1)m}\chi(m) e^{[x+m]_q t}-t
\sum_{m=0}^\infty q^{hm +x} \chi(m) e^{[x+m]_q t }. \eqno(13) $$
Therefore, we obtain  the  following generating function:

\medskip
{ \bf Proposition 9.} Let $      F_{\chi, q}(t, x  \mid h)  =
\sum_{n=0}^\infty \beta_{n, \chi,
 q}^h(x)
\dfrac{t^n}{n!}.$  Then we have
$$   F_{\chi, q}(t, x  \mid h) =(h-1)(1-q)
\sum_{m=0}^\infty q^{(h-1)m}\chi(m) e^{[x+m]_q t}-t
\sum_{m=0}^\infty q^{hm +x} \chi(m) e^{[x+m]_q t }.  $$
\medskip

 Let $     F_{\chi, q}(t \mid h)  = \sum_{n=0}^\infty \beta_{n, \chi,
 q}^h\dfrac{t^n}{n!}.$  Then we also get
$$   F_{\chi, q}(t \mid h) =(h-1)(1-q)
\sum_{m=0}^\infty q^{(h-1)m}\chi(m) e^{[m]_q t}-t
\sum_{m=0}^\infty q^{hm } \chi(m) e^{[m]_q t }. \eqno(14) $$ From
(14) and the definition of $ \beta_{n, \chi,
 q}^h$, we can derive the
  following functional equation:
$$ \aligned  \beta_{n, \chi,
 q}^h & = \left.  \frac{d^n}{dt^n} F_{ \chi, q}(t \mid h)
 \right\vert_{t=0}\\
& =  (h-1)(1-q) \sum_{m=0}^\infty q^{(h-1)m} \chi(m)[m]_q^n -n
\sum_{m=0}^\infty q^{hm} \chi(m) [m]_q^{n-1}. \endaligned
\eqno(15)$$ By (14) and (15), we obtain the following corollary:

\bigskip
 { \bf Corollary 10.} For $ h \in \Bbb C, n \in \Bbb Z_+$ ,   we have
$$   \beta_{n, \chi,
 q}^h =  (h-1)(1-q) \sum_{m=0}^\infty q^{(h-1)m} \chi(m)[m]_q^n -n
\sum_{m=0}^\infty q^{hm} \chi(m) [m]_q^{n-1}.
$$
\bigskip

For  $ s \in \Bbb C \setminus \{ 1 \}, h \in \Bbb C, $ and  $ x
\neq 0, -1, -2, \cdots, $   we consider  complex  integral as
follows:
$$ \aligned   & \dfrac{1}{\Gamma(s)} \int_0^\infty  F_{\chi, q}(-t, x \mid h)t^{s-2}dt \\
&=  \dfrac{(h-1)(1-q)}{s-1}  \sum_{m=0}^\infty \dfrac{ q^{(h-1)m
}\chi(m)}{ [x+m]_q^{s-1}} + \sum_{m=0}^\infty \dfrac{ q^{hm +x}
\chi(m)}{ [x+m]_q^s}.
\endaligned
\eqno(16)
 $$
From (16), we can define Hurwitz's type $(h, q)$-$L$-function as
follows:

\bigskip
{ \bf Definition 11.}
 For  $ s \in \Bbb C \setminus \{ 1 \}, h \in \Bbb C, $ and  $ x
\neq 0, -1, -2, \cdots, $  define
$$ L_q^h( s, \chi \mid x)=  \dfrac{(h-1)(1-q)}{s-1}  \sum_{m=0}^\infty \dfrac{ q^{(h-1)m }\chi(m)}{
[x+m]_q^{s-1}}+ \sum_{m=0}^\infty \dfrac{ q^{hm +x}\chi(m)}{
[x+m]_q^s}.
 $$
\bigskip

By the definition of the generating function for the generalized
$(h, q)$-Bernoulli polynomials attached to $\chi$, we get
$$ \dfrac{1}{\Gamma(s)} \int_0^\infty  F_{\chi, q} (-t, x \mid h)t^{s-2}dt
= \sum_{n=0}^\infty \dfrac{(-1)^n \beta_{n, \chi, q}^h(x)}{n!}
\dfrac{1}{\Gamma(s)} \int_0^\infty t^{n+s-2}dt . \eqno(17)$$
 We
see form (16) and (17) and the  basic theory of complex analysis
including Laurent series that
$$
\aligned & -\dfrac{(h-1)(1-q)}{k}  \sum_{m=0}^\infty   q^{(h-1)m }
\chi(m) [x+m]_q^k + \sum_{m=0}^\infty  q^{hm +x} \chi(m)
[x+m]_q^{k-1} \\
&  = - \dfrac{ \beta_{k, \chi, q}^h(x)}{k}, \text{ for } k \in
\mathbb{N }. \endaligned \eqno(18)$$ From (18) and Definition 2,
we obtain the following functional equation.

\medskip
{ \bf Theorem 12.}
 For $ k  \in \Bbb N, $ we
have
$$ L_q^h( 1-k, \chi \mid x)=  - \dfrac{ \beta_{k, \chi,
q}^h(x)}{k}.
 $$
\medskip

Let $\chi$ be non-trivial   Dirichlet  character with conductor $f
\in \Bbb N$. Then  we can also consider Dirichlet's type
$(h,q)$-$L$-function  as follows:
$$ L_q^h( s, \chi)=  \dfrac{(h-1)(1-q)}{s-1}  \sum_{m=1}^\infty \dfrac{ q^{(h-1)m }\chi(m)}{
[m]_q^{s-1}}+ \sum_{m=1}^\infty \dfrac{ q^{hm }\chi(m)}{ [m]_q^s},
\eqno(19) $$ for $ s \in \Bbb C \setminus \{ 1 \}, h \in \Bbb C.$

By Corollary 10 and (19), we get
$$ L_q^h( 1-k, \chi)=  - \dfrac{ \beta_{k, \chi, q}^h}{k}, \text{ for } k \in
\mathbb{N }. $$

\end{document}